\documentclass{amsart}
\usepackage[pdftex,colorlinks]{hyperref}

\newcommand{\y}{\textsf{yes}}
\newcommand{\n}{\textsf{no}}

\newcommand{\PP}{\textsf{P}} 
\newcommand{\RR}{\mathbb{R}}

\newcommand{\NN}{\mathbb{N}}

\begin{document}

\title{On foundations for deductive mathematics}
\author{Frank Quinn}

\date{June 2024}

\maketitle


\section{abstract}\label{sect:abstract} 
This article was motivated by the discovery of a  potential new foundation for mainstream mathematics. The goals are to clarify the relationships between primitives, foundations, and deductive practice; to understand how to determine what is, or isn't, a foundation; and get clues as to how a foundation can be optimized for effective human use.
 For this we turn to history and professional practice of the subject. We have no asperations to Philosophy.  

The first section gives a short abstract discussion, focusing on the significance of consistency. 
The next briefly describes  foundations, explicit and implicit, at a few key periods in mathematical history. We see, for example, that at the primitive level human intuitions are essential, but can be problematic. We also see that traditional axiomatic set theories, Zermillo-Fraenkel-Choice (ZFC) in particular, are not quite consistent with mainstream practice. 
The final section sketches the proposed new foundation and gives the basic argument that it is uniquely qualified to be considered {the} foundation of mainstream deductive mathematics. The ``coherent limit axiom'' characterizes the new theory among ZFC-like theories. This axiom plays a role in recursion, but is implicitly assumed in mainstream work so does not provide new leverage there. In principle it should settle set-theory questions such as the continuum hypothesis. 

\section{Abstractions} ``Foundation'' refers to a body of deductive practice, here mainstream deductive mathematics. Historically practice comes first. Foundations are a way to analyze the body of practice, particularly its consistency. 
\begin{enumerate}\item A \textbf{primitive} is an object,  hypotheses, or method of argument that is postulated rather than deduced from something else.
 \item A primitive \textbf{foundation}  is a  collection of primitives that is complete in the sense that no other primitives are needed.  Note that the methods of argument (the logic) is part of the foundation. 
 \item  A foundation is \textbf{consistent}  if every chain of deductions  that produces a contradiction can be shown to have a logical error.
\end{enumerate}

\subsection{About consistency}\label{ssect:consistency}
 The operational view of consistency is ``complete reliability'', in the sense that  the outcome of an error-free argument will never be contradicted by the outcome of any other such argument. Moreover, the outcome will not introduce errors if used as an ingredient in further arguments. Arguments by contradiction for example,  assume something,  deduce a contradiction, and conclude that the thing assumed must be false. Success is sensitively dependent on \emph{everything else} used in the argument being completely reliable, and that the logic itself does not introduce errors. 
 
Consistency provides an internal criterion for correctness. In practice, however, it is far more important as a criterion for \emph{in}correctness. In reasonably ambitious research programs, 95 percent or more of attempted proofs turn out to have flaws. Sometimes, with luck and persistence, these flaws can  be corrected, and the attempts converge to a valid argument. Sometimes years of effort have to be abandoned as hopeless, or outright wrong. The payoff for this painful process is closure: even outrageously counterintuitive conclusions are accepted when the arguments are carefully checked and found to be error-free. 

\subsubsection*{Empirical} G\"odel has shown that  foundations generally cannot prove their own consistency. However, we can \emph{test} consistency: make lots of deductions, and verify that any contradictions are accounted for by logical errors.  Sometimes inconsistencies do not appear until work reaches a certain level of sophistication. Another wrinkle is that a foundation can be technically inconsistent but consistent \emph{in practice} in the sense that methodology has evolved so as to steer people away from the inconsistencies. 
The main point is that empirical consistency is  a property of a methodology rather than of a foundation. 
 Explicit foundations are largely retrofits that record successful experience and narrow down potential sources of error.   
 
\subsubsection*{Global} Evidence for consistency applies to the development as a whole, not just individual results.  In other sciences, and in mathematics when trying out new methodologies, experiment can give evidence that individual conclusions are reliable, within limits and independently of most other conclusions.  Foundational consistency checks are global in that they support consistency of the whole development. This is the flip side of the fact that a single genuine failure would throw doubt on the whole development. 

\subsection{Definitions and axioms vs  primitives} We clarify  usage in the mathematical community.
 
 \textbf{Definitions} are located inside a developed context, and usually consist of \textbf{axioms} that the things being defined should satisfy.  Properties of these things are to be inferred from the  axioms and more-basic ingredients. In particular, while intuition and intent can guide the formation of a definition they
 have no logical force. Indeed, first attempts at definitions typically fail to have intended properties. Major standard definitions are often the result of years of work to fix failures and consolidate partial successes, and should be thought of as distilled hard-won wisdom. They are certainly not arbitrary constructs. 
 
\textbf{Primitive} objects, hypotheses, etc.~are starting points for the logical context, \emph{not} defined in terms of lower-level objects. Primitive hypotheses, in particular, are not axioms in the contemporary sense. Since properties cannot be inferred from lower-level material, they depend on heuristic ideas, syntax specifications, and examples. 

For example,  Euclid's description of points as ``things without width or length'' is not a definition in the contemporary sense because ``width'' and ``length'' have not been given precise meaning. It succeeds as a primitive because this, filtered through physical experience, some pictures, and examples of usage, reliably evokes consistent understandings. Polygons are `defined' because they are described in terms of primitives. 

We emphasize  that consistency at the primitive level depends heavily on different users extracting functionally equivalent understandings. This has often been a problem. 
\subsection{Deeper is better} In practice, definitions are often easier to learn and use with precision than are complicated primitives. For example,  attempts to use the real numbers as a primitive were problematic. Now the reals are defined in terms of natural numbers, and natural numbers are defined using set theory. This approach has made full-precision understanding of the reals faster and more routine.  

Curiously,  another benefit of working several levels above the primitives is that it can filter out misunderstandings or ambiguities about the primitives. For instance, the traditional axiomatic approach to set theory takes `sets' as primitives, and the primitive hypotheses proscribing their behavior are quite complicated. But these do not have to be understood in detail to get a completely solid description of the real numbers. The new approach to set theory has a similar advantage: the primitives (object generators) are several levels below sets, so sets are \emph{defined} and their properties proved rather than hypothesized. This gives a more accessible, and slightly sharper, understanding of set theory. 
\subsection{Other views} 
There is a vast literature, mostly Philosophical, about the nature and role of foundations.  It is far from consistent, and while there are interesting insights, cf.~\cite{maddy}, 
 some of it is pretty silly. It does not seem to be useful to try to sort through it here.

\section{Sketchy history}\label{sect:history} We draw lessons from  a few key historical deveopments, mainly to provide perspective on what happened in the early twentieth century. The accounts are ideosyncratic, brief, and superficial, but sufficient for the purposes here. For a bit more detail I have found the  essays in `The Stanford Encyclopedia of Philosophy'  to be helpful, for example \cite{bagara} and \cite{ferreiros}. 

\subsection{Ancient Greece}
The first relatively explicit use of a foundation was in  the geometry and number theory of ancient Greece. 
The primary significance at the time seems to have been its use as a metaphor for order in the natural world, not for its ``real-world'' applications.  Everyone could see, and be impressed by, genuine consistency coming from logic and innate understanding, as opposed to the chaotic array of beliefs due to cultural bias or religious or philosophical doctrine. The practical applications are pretty weak. To some extent this can be seen as a quantifier issue: for metaphors it is sufficient that the theory have \emph{some} striking  outcomes. Practical use requires ways to describe outcomes for \emph{all} problems in appropriate settings. For example, clever geometric tricks can determine the area between a line and an arc of a circle in a few special cases. Calculus gives a description of \emph{every} such area, and the appearence of transcendental functions in the solution shows that very few cases can be solved with ancient techniques. Nonetheless many people find the classical tricks more meaningful than  routine but effective calculus.

 \subsubsection*{Meaningful vs effective} We expand on this last remark.
It has been known for over fifty years that humans have an innate physics, and it is essentially Aristotlean. People therefore find
Aristotlean physics  comfortable and meaningful, and it was accepted for centuries even though it fails very simple consistency tests. In general, humans  seem not to expect consistency, either internal or with reality, of the things we find meaningful. 

The lesson for geometry is that  Euclidean geometry is in large part an exploration of human visual perception and physical experience. The fact that it accurately reflects the structure of the plane is because our visual perception does, not because it is a direct exploration of the plane. Much of its power as a metaphor comes from this coincidence: it is comfortable because it matches our perceptions, and effective because in this case perception accurately reflects reality.  Another lesson from Aristotlean physics is that outside of some areas in geometry our physical intuitions do not match reality very well. Approaching mathematics through direct perception seldom leads to effective internalization, but is attractive and comfortable enough that it remains a strong theme in elementary education.
\subsection{Europe in the 1600s} 
The next major shift began in the early 1600s when Galileo, Kepler, and others showed that mathematics could be directly powerful in the world, not just a metaphor. Mathematical development became driven by the needs of physical models, and methodology evolved rapidly. For instance, there is a conflict between ratios and negative numbers that had inhibited use of negatives. This was resolved by replacing ratios with fractions. Compact algebraic notation replaced earlier discursive formulations. The ancient Greeks were justly proud of their resolution of ``multitudes'' into natural numbers, and units indicating the things being counted, but the analogous resolution of ``magnitudes'' into real numbers and units was finally accomplished only in the 1600s. 
 By the time Newton, Leibnitz, and their contemporaries were active, the power of application-oriented mathematics was firmly established and core methodologies were taking shape.

Infinitesimals were often part of the foundational intuitions of the real numbers in this period. They worked quite well for a time, but eventually became problematic and had to be replaced by limits. We explain how this fits into the analysis here. Consider the truncated polynomial ring \(\RR[\delta]/(\delta^2=0)\). The implicit presumption was that a function \(f\) of a real variable should extend to this ring, and \(f[x+y\delta]=f[x]+g[x]y\delta
\) where \(g\) is the derivative of \(f\). This is true if \(f\) is given by a power series: simply plug the extended variables into the series. Infinitesimals thus worked as long as it was harmless to assume functions were piecewise analytic.  But limits of analytic functions need not even be piecewise continuous, and continuous limits can be nowhere differentiable. When limits became important (eg.~for Fourier and Laplace series) infinitesimals became untenable. The methodology had to evolve from one in which differentiability was assumed, to a more primitive approach in which  differentiability had to be defined and proved. 

In the 1960s Abraham Robinson  showed infinitesimals can be used  consistently if  arguments are restricted to first-order logic. But this first-order logic constraint is strongly inconsistent with standard practice. This is a potential special-purpose tool, not an alternative general approach. 

\subsection{Mid to late 1800s} In this period it was still generally felt that real numbers should be understood through intuitions from physical experience, though, of course, this should not include earlier intuitions about infinitesimals.  This became problematic as goals became more ambitious. People without sufficiently accurate intuitions had unreliable outcomes, and were reduced to appeals to authority (``Gauss did it so it must be ok'') or were shut out of the game. From a foundational point of view, ``intuitive reals'' became unsatisfactory as a primitive because it was not uniformly understood by different users. This is a bit different from problems with infinitesimals, which in the larger context are inconsistent even when correctly understood. 

 In the mid 1800s Dedekind  rigorously  described the reals in terms of natural numbers. This made them reliably accessible but didn't fully fix the difficulties,  because the natural numbers had varying interpretations. Some  intuitively-attractive ones  lacked logical force. This was addressed by defining natural numbers using early versions of set theory. For human use, sets seemed to be effective starting points. By 1880 Dederkind was using set-theoretic methods and terminology, and the idea that mathematics could have a foundation in a theory of sets began to have traction;  cf.~\cite{ferreiros}. 
 
 By 1900, foundation-oriented mathematicians were having good success with what we now know as na\"\i ve set theory. 
Coming from another direction, Frege, Cantor, and others clarified intuitive notions about set theory, and showed that it had its own substance.
\subsection{Early 1900s} 
In 1902 Russell publicized his paradox showing na\"\i ve set theory to be inconsistent. Reactions were quite different in the mainstream and set-theory communities, and resulted in them drifting apart.
\subsubsection*{Mainstream community}
 By 1900 quite a few mathematicians had found na\"\i ve set theory to be effective and reliable, and many leading mathematicians were satisfied that all mainstream work could be derived from it.   The response to Russell's paradox was to add the caution ``don't say `set of all sets' ''. This did not effect actual practice since nobody saw a need to do that anyway. 
 
By the 1920s much of the younger generation had embraced the set-theory foundation  and the rigor it enables, and were extremely successful with it.  The main set-theory concern of the period was that the axiom of Choice seemed too good to be true, not that there would be more problems  with na\"\i ve set theory.

\subsubsection*{Set-theory community}
The paradox seriously set back the study of sets as a subject. After a few more unsuccessful attempts they gave up trying to define sets directly, and instead began abstracting key properties (of na\"\i ve set theory) to use as axioms.

In the next phase of the development, the notion of ``function'' seemed problematic. A na\"\i ve function \(X\to Y\) is ``an assignment of an element of \(Y\) to each element of \(X\)''. This does not qualify as a definition because ``assignment'' is neither defined nor declared as a primitive. In fact this description of `function' works well as a primitive, but the community choose not to follow this route. Instead a single function, the `membership' operator was taken as a primitive. Specifically, an implementation (or ``model'') of the axioms is a pair \((U,\in)\), where \(U\) is a `universe'  of possible elements, and the membership operator is a binary pairing \(\in\colon U\times U\to \y/\n\). In general, subcollections of \(U\) correspond to binary functions \(U\to \y/\n\). The \emph{sets} of the theory correspond to binary functions of the form \((\#\in x)\) for some \(x\in U\). Note that elements are used to parameterize sets. One source of confusion is that this parameterization is sometimes described as an identification: ``elements of sets are again sets''. This phrasing may be a hangover from 19th century philosophy. More general functions  are then defined in terms of sets. Specifically, a function \(X\to Y\) is a correspondence whose graph is a sub\emph{set} of \(X\times Y\). The standard forms of the axioms up to this point are mainly due to Zermillo. 

A problem with the Zermillo axioms is that they do not ensure there are enough functions to transact the basic business of set theory. Experimentation, drawing on a strong legacy of formal logic, revealed that  functions obtained from set operations by first-order logic are sufficient to give workable rigorous set theories. Standard formulations of the logical hypotheses are mainly due to Fraenkel. Many variations have been considered, but after about 1920  the Zermillo-Fraenkel (ZF) axioms were widely accepted as standard.    We will see, however,  that while this  solved a problem in set theory, it disconnected formal set theory from the mainstream. 
 
\subsection{Mid 1900s to early 2000s} 
Early in this period the axiom of Choice was accepted as well-tested, and included in standard formulations of both communities. Zermillo-Fraenkel-Choice became the ``gold standard'' in the set-theory community, and na\"\i ve-with-choice (as always, with the set-of-sets constraint) remained the implicit standard for the mainstream. 
\subsubsection*{Mainstream community}
 Na\"\i ve set theory---as used in standard practice---enabled highly rigorous arguments of a depth unimaginable in the previous century. The consistency of na\"\i ve set theory---again, as used in standard practice---is extremely well tested. It is possibly our most solidly established empirical conclusion. There are some drawbacks, but there is no reason to accept any loss of functionality to address these drawbacks. 
 
One drawback is that  category theory does not fit comfortably in standard set theories.  Categories were introduced essentially as bookkeeping devices to encode large-scale structures in algebraic topology and homological algebra. For these purposes ignoring set-theory mismatches seemed harmless.  By the late 1900s delicate work with higher-order category theory and point-set properties of universal spaces needed more precision. There were proposals to use categories as a substitute for set theory, cf.~\cite{marquis}, \cite{maddy}. Grothendieck developed a variation, topos theory, initially as a setting for vigorous abstract algebraic geometry. As it matured this was also promoted as a substitute for set theory. 
Later Voevodsky proposed using his ``univalent foundation'' as a setting for mathematics, citing its effectiveness for computer implementation. 

These proposals reflect a misapprehension about foundations. We have stressed that the core job of a foundation is to give a distilled starting point that identifies key consistency checks, and to encode the results of such checks. To be effective a foundation should be  minimal and as simple as possible. The proposals above are quite elaborate, and designed to display mathematical structure or make it more mechanically accessible. These are good things for applications but, as we saw with the real numbers, working directly with high-level concepts is likely to cause problems. The situation seems to be that for generations,  set theory has provided such a rock-solid base for mainstream practice  that we have forgotten what it is like to have to worry about the reliability of our tools.  The proposals   implicitly assume consistency will follow from our experience with set theory. Indeed many of them implicitly depend on  a background set theory, so implicitly presume consistency to be inherited from that. 

The foundation described in the next section begins several levels below set theory. The lower-level environment seems to give an effective setting for categories. 
\subsubsection*{Set-theory community}
The ZFC axioms were found to be under-determined, and in 1963 Paul Cohen introduced Forcing as a systematic way to construct new models. A flowering of the subject ensued, with many thousands of pages devoted to showing a wide variety of hypotheses are consistent with, or independent of, the basic axioms. Deep results were found in other directions, PCF theory and large cardinal axioms to mention only two. See \cite{jech} for an extensive introduction. 

However, these theories  do not  provide a foundation for mainstream practice. Sufficiency of na\"\i ve set theory was asserted  around 1900, and this assertion does \emph{not} apply to  theories that require first-order logic. The general problem is described by the ``coherent limit axiom of \cite{quinn1} but we illustrate it in an important special case.Recall that every real number in the unit interval has a base-2 decimal expansion that assigns either 0 or 1 to each natural number. This assignment can also be thought of as specifying a subcollection of \(\NN\). The restriction of the assignment to any finite interval \([0,n]\) of natural numbers is a function in any ZFC theory. Equivalently, the intersection of the subcollection with \([0,n]\) is a subset in any ZFC theory. The full assignment is the limit of these restrictions, or equivalently the subcollection is the union of the intersections. However, there are ZFC models in which some of these limits are not functions in the model, or equivalently the union of the finite subsets is not a set in the model. The corresponding real numbers are missing from such theories. To insure a particular real number is in a particular model we would have to employ either special information about the model, or  first-order logic. This is an issue even for elementary calculus. But in the mainstream nobody addresses this. The conclusion is that ZFC is too restrictive to describe standard mathematical practice, even at basic levels. 

 The insufficiency of ZFC  does not seem to be a deep point. The fact that it went generally unremarked may be another instance of meaningfulness overshadowing concern for effectiveness.  

\section{Object generators} We sketch the proposed new foundation;  see \cite{quinn1} for details, and \cite{quinn2} for a guide to routine use. 

\subsection{Generators} The full story begins several levels below set theory, with ``object generators'' and their morphisms.
The intuition is that object generators ``generate objects'', essentially like the `object' primitive in category theory. If \(G\) is an object generator then \(X\in\in G\) means ``\(X\) is a an object in \(G\)'', and generators are defined using the syntax ``\(X\in\in G\) means (\(\cdots\))''. 

Essentially nothing else is included. Unlike sets, generators provide no way of identifying their own outputs, and no way of knowing if two outputs are the same. 
Further, while expected uses are for mathematical objects, for instance \(X\in\in\textbf{ groups }\) means ``\(X\) is a set with a group structure'', nothing prevents silly examples, for instance \(X\in\in G\) means ``\(X\) is a Tuesday in the month of May, 1917''. Rather than try to prevent such things we filter them out at a later stage. 
Note that few requirements means few opportunities for contradiction. 

``Morphisms'' of generators are essentially the same as ``functions'' in na\"\i ve set theory: a morphism \(F\colon G\to H\) is an assignment of an object \(F[X] \in\in H\) to every object \(X\in\in G\). 

\subsection{Logic} 
The native logic of object generators is   non-binary in the sense that we might \emph{assert} that \(X,Y\) are the same, but in general there is no (yes, no)-valued function that can detect this. This logic is unfamiliar and somewhat complicated, and the first step toward set theory is to establish a sub-context that does use binary logic.

The key ingredients  are \textbf{binary functions}: functions to an object generator with exactly two objects. We denote this by  \(\{\y, \n\}\), or \(\{0,1\}\). One of the primitive hypotheses (axioms) asserts the existence of  a two-object generator. Standard binary logic applies to such functions. Incidently, these functions have no allowance for time dependence, so they cannot detect things in the physical world. This is how silly examples of object generators get filtered out. 
\subsection{Domains and sets}

A \textbf{logical domain} is an object generator, say \(A\), with a binary pairing \(A\times A\to\{\y,\n\}\) that returns `yes' if the two elements are the same. The name `logical domain' is supposed to suggest that these are natural settings for binary logic. 

Quantification is the final ingredient. If \(A\) is a logical domain then the object generator whose objects are binary functions on \(A\) is called the ``powerset'' of \(A\), and is denoted by \(\PP[A]\). The domain \textbf{supports quantification} if there is a binary function \(\PP[A]\to \{\y,\n\}\) that returns `yes' if and only if the input is the empty (always-`no') function. In standard logic the empty-detecting operator is given by \(f\mapsto \forall x\in A, (f[x]=\n)\). The point is that if this one quantification expression is implemented by a binary function, then all quantifications over \(A\) are. 

Finally, \textbf{sets} (\emph{relaxed} sets, if distinctions are necessary) are defined to be logical domains that support quantification. 
 \subsection{Hypotheses} There are four primitive hypotheses. Three are standard: there is a 2-element set; the axiom of Choice; and the natural numbers support quantification. To clarify this last, the first two hypotheses suffice to construct the natural numbers, but do not imply that it supports quantification. This is equivalent to the standard Axiom of Infinity in ZFC. The final hypothesis asserts that if a logical domain supports quantification then so does its powerset. 
 
Note that ZFC has many more, and more complex, axioms.
\subsection{Comparisons}
We compare this theory, na\"\i ve set theory,  ZFC axiomatic set theories, and standard mathematical practice. First, the ``='' operator built into the logic of na\"\i ve and ZFC allows comparison of any two elements of any two sets. In relaxed set theory ``='' is only defined for elements of a single set. Mathematical practice is modeled on na\"\i ve set theory so in principle a global ``='' is available but, as far as the author can tell, it is never used.

Functions are primitives in relaxed and na\"\i ve set theory; are described in the same way; and both are consistent with standard practice. In ZFC only the equality and membership functions are primitive. Functions obtained by first-order logical expressions are also functions in the theory. Beyond this, functions are defined to be correspondences whose graphs are subsets. As explained in the previous section, this  is not consistent with standard practice. 

Finally, relaxed sets satisfy versions of the ZFC axioms. This is made precise in \cite{quinn1} with the construction of an object that satisfies the ZFC-1 axioms, where ``-1'' means ``ignore first-order logical constraints''. As above, this makes sense because functions are taken as primitives. This construction turns out to be  maximal so, in particular, all models of ZFC uniquely embed in this as transitive sub-models. The ``coherent limit axiom'', a generalization of the real-number problem described at the end of the last section, is shown in \cite{quinn1} to characterize this maximal theory. This axiom is routinely assumed in standard practice so standard practice, by default, takes place in relaxed set theory. 
We note that identifying the `missing axiom' (G\"odel's phrase) does not lead to new mainstream methodology, but rather emphasizes the extra work that would be required to stay in a smaller ZFC model.

The conclusion is that the relaxed set theory defined in \cite{quinn1} is uniquely qualified to be a foundation consistent with mainstream mathematical practice.


\begin{thebibliography}{quinn 2}

\bibitem[Bagaria]{bagara}Joan Bagaria, \emph{Set Theory}, The Stanford Encyclopedia of Philosophy (Winter 2023 Edition), Edward N. Zalta and Uri Nodelman (eds.), URL = https://plato.stanford.edu/archives/win2023/entries/set-theory/ .


\bibitem[Ferreir\' os]{ferreiros} Jos\'e Ferreir\'os, \emph{The Early Development of Set Theory},  The Stanford Encyclopedia of Philosophy (Summer 2023 Edition), Edward N. Zalta and Uri Nodelman (eds.), URL = https://plato.stanford.edu/archives/sum2023/entries/settheory-early/ . 

\bibitem[Marquis]{marquis}Jean-Pierre Marquis, "Category Theory", The Stanford Encyclopedia of Philosophy (Fall 2023 Edition), Edward N. Zalta and Uri Nodelman (eds.), URL = https://plato.stanford.edu/archives/fall2023/entries/category-theory/ . 

\bibitem[Jech]{jech} Thomas J.~Jech:
\emph{Set theory -- 3rd Millennium ed, revised and expanded.}
Springer monographs in mathematics (2002)
ISBN 3-540-44085-2

\bibitem[Maddy]{maddy} Penelope Maddy \emph{What do we want a foundation to do?
Comparing set-theoretic, category-theoretic, and univalent approaches} Springer nature 2019,
ebook ISBN 978-3-030-15655-8

\bibitem[Quinn 1]{quinn1} Frank Quinn, \emph{Object generators,  relaxed sets, and a foundation for mathematics} preprint 2023, http://arXiv.org2110.01489

\bibitem[Quinn 2]{quinn2} Frank Quinn, \emph{Object generators,  categories, and  everyday set theory} preprint 2023, http://arxiv.org/abs/2210.07280

\end{thebibliography}
\end{document}